\documentclass[graybox]{svmult}

\usepackage{mathptmx}        
\usepackage{helvet}          
\usepackage{courier}         
\usepackage{makeidx}         
\usepackage{graphicx}        
\usepackage{multicol}        
\usepackage[bottom]{footmisc}

\makeindex             

\usepackage{amsmath}
\usepackage{amssymb}

\usepackage{amsthm}
\usepackage{relsize}
\usepackage{bm}
\usepackage{booktabs}
\usepackage{tabularx}

\theoremstyle{plain}
\newtheorem{Th}{Theorem}[section]

\newtheorem{Prob}[Th]{Problem}

\theoremstyle{remark}
\newtheorem*{Rem}{Remark}

\newcommand{\R}{\mathbb{R}}
\newcommand{\Tensor}{\R^{3\times 3}}
\newcommand{\SymTensor}{\R^{3\times 3}_{\mathsmaller{\mathrm{sym}}}}
\newcommand{\uT}{^{\mathsmaller{\mathrm{T}}}}
\newcommand{\dd}{\mathop{}\!\mathrm{d}}
\newcommand{\vep}{\varepsilon}
\newcommand{\GammaD}{\Gamma_{\mathsmaller{\mathrm{D}}}}
\newcommand{\GammaN}{\Gamma_{\mathsmaller{\mathrm{N}}}}
\newcommand{\fN}{f^{\mathsmaller{\mathrm{N}}}}
\newcommand{\Jump}[1]{\langle #1 \rangle_S}

\def\1/2{\frac{1}{2}}

\DeclareMathOperator{\Div}{div}

\begin{document}

\title*{An energy-consistent model of dislocation dynamics in an elastic body}
\author{Vladim\'ir Chalupeck\'y and Masato Kimura}
\institute{%
Vladim\'ir Chalupeck\'y \at Fujitsu Limited, 1-17-25 Shinkamata, Ota-ku, Tokyo 144-8588, Japan\\%
\email{chalupecky@jp.fujitsu.com, vladimir.chalupecky@gmail.com}%
\and%
Masato Kimura \at Faculty of Mathematics and Physics, Kanazawa University, Kanazawa 920-1192, Japan\\%
\email{mkimura@se.kanazawa-u.ac.jp}%
}
\maketitle

%

\abstract*{We propose an energy-consistent mathematical model for motion of
dislocation curves in elastic materials using the idea of phase field model.
This reveals a hidden gradient flow structure in the dislocation dynamics.  The
model is derived as a gradient flow for the sum of a regularized Allen-Cahn
type energy in the slip plane and an elastic energy in the elastic body. The
obtained model becomes a 3D-2D bulk-surface system and naturally includes the
Peach-Koehler force term and the notion of dislocation core.  We also derive a
2D-1D bulk-surface system for a straight screw dislocation and give some
numerical examples for it.}

\abstract{We propose an energy-consistent mathematical model for motion of
dislocation curves in elastic materials using the idea of phase field model.
This reveals a hidden gradient flow structure in the dislocation dynamics.  The
model is derived as a gradient flow for the sum of a regularized Allen-Cahn
type energy in the slip plane and an elastic energy in the elastic body. The
obtained model becomes a 3D-2D bulk-surface system and naturally includes the
Peach-Koehler force term and the notion of dislocation core.  We also derive a
2D-1D bulk-surface system for a straight screw dislocation and give some
numerical examples for it.}

\section{Introduction}\label{sec:1}

A dislocation or a dislocation curve is a crystallographic line defect within a
crystal structure which was first studied by E.~Orowan, M.~Polanyi and
G.I.~Taylor independently in 1934. It is considered as a main mechanism of
plastic deformation or yielding of the material and various material properties
are studied in relation to the dislocation dynamics nowadays. See
\cite{Bulatov06, Nabarro67} and references therein for more details.

Some properties of dislocations are as follows. In many cases, a dislocation
is a plane curve in a fixed slip plane, and is the boundary of a region shifted
by the Burgers vector $b\in\R^3\setminus\{0\}$. The Burgers vector is tangential
to the slip plane and coincides with a translation vector of the crystal
structure. Although end points often appear on the material surface or at other
defects inside in a real crystal material, it is known from a topological
argument that a dislocation curve cannot have an end point inside of the
crystal lattice. Its typical length is around $10^{-6}$ m, typical thickness is
around $10^{-9}$ m, and typical velocity is $10^{-6}$--$10^{2}$ m/s.

A dislocation generates a stress field around it by deforming the crystal
structure, and interacts with other dislocations, defects, and far field
conditions through the stress field. The virtual force from the stress field
acting on the dislocation is called the Peach-Koehler force \cite{Bulatov06}.

There are some mathematical models of dislocations (see \cite{Bulatov06} and
references therein), however, the following mathematical difficulties exist.
The displacement $u$ has a jump along the dislocation curve and it has infinite
elastic energy ($u\not\in H^1_{\mathrm{loc}}$). Therefore the elasticity equation is
valid only outside of the dislocation core of radius about $5|b|$, where $|b|$
corresponds to the interatomic spacing.

In this paper, we construct an energy-consistent model in a mathematically
clear way by
means of the idea of the phase field \cite{Carter97}. We suppose a
quasi-stationary condition, which means that the small deformation and the
stress field of the material are described by the static linear elasticity
equations. In other words, they are given by a minimum energy state of a
suitable elastic energy. Our phase field model for the dynamics of a
dislocation will be derived as a gradient flow of a total energy including the
elastic one, as an analogy to the Allen-Cahn equation \cite{Chen92,Fife88}.

\begin{figure}[b]
\begin{center}
\includegraphics[width=0.7\linewidth]{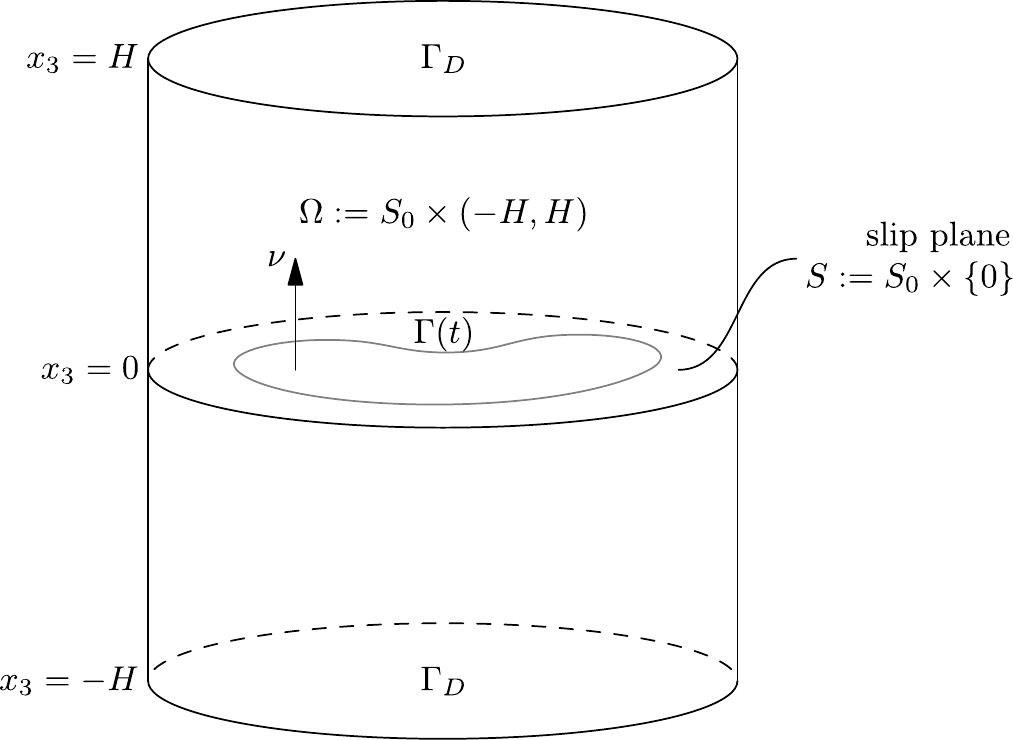}
\end{center}
\caption{A typical configuration of an elastic body with a closed dislocation
curve in a slip plane.}
\label{fig:1}
\end{figure}

\begin{figure}
\sidecaption
\includegraphics[width=0.4\linewidth]{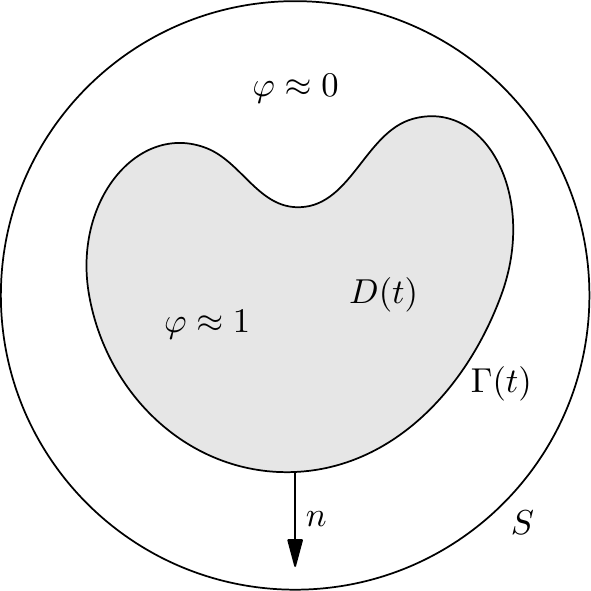}
\caption{A slip plane with a dislocation curve $\Gamma(t)$. The phase field
$\varphi$ approximately satisfies $\varphi\approx 1$ in $D(t)$ and
$\varphi\approx 0$ in $S\setminus D(t)$.}
\label{fig:2}
\end{figure}

We denote a dislocation curve by $\Gamma(t)$. In this paper, we assume the
dislocation curve $\Gamma(t)$ is a smooth Jordan curve in a plane. The plane
which includes the dislocation curve is called the slip plane. The crystal
lattice has a defect along the dislocation curve and the defect is represented
by the Burgers vector $b\in\R^3\setminus \{0\}$. We suppose that $b$ is
tangential to the slip plane. Let $\tau\in\R^3$ be a counterclockwise
tangential unit vector of the dislocation curve $\Gamma(t)$ and let $n\in \R^3$
be one of the two possible directions of the unit normal vector 
of $\Gamma(t)$ tangential to the slip plane.
We call the directions $n$ and $-n$ ``outward'' and ``inward'', respectively,
in this paper.
We choose a unit normal vector $\nu\in\R^3$ of the slip plane with which
$(n,\tau,\nu)$ is a right-handed orthonormal coordinate system.  See
Fig.~\ref{fig:1} and Fig.~\ref{fig:2} for a typical configuration.

One of the simplest mathematical models for the motion of a dislocation curve is
\cite{Minarik10}
\begin{equation}\label{cf}
    V = -c\kappa + f, \qquad \text{on} \quad\Gamma(t),
\end{equation}
where $V$ is the outward normal velocity of $\Gamma(t)$ in the slip plane and
$\kappa$ is the inward signed curvature of 
the dislocation when considered as a plane
curve. The equation \eqref{cf} is called the mean curvature flow and is a typical
mathematical model of phase transition. It is also known that it appears as a
singular limit problem of the Allen-Cahn equation \cite{Chen92,Fife88}.

The term $f$ in \eqref{cf} represents a far-field interaction through an
elastic field from other dislocation curves including $\Gamma(t)$ itself, other
defects, and boundary conditions. It is known that there is a virtual force
$F\in\R^3$ acting on the dislocation curve which is given by the Peach-Koehler
formula \cite{Bulatov06}
\begin{equation}\label{PK}
    F = \tau\times (\sigma b),
\end{equation}
where $\sigma\in\SymTensor$ is the stress tensor field, $b\in\R^3$ is the Burgers
vector, and $\tau\in \R^3$ is the unit tangential vector of $\Gamma(t)$. Then the
term $f$ in \eqref{cf} is formally given by the $n$-direction component of the
Peach-Koehler force $F$. Hence, we obtain
\begin{equation}\label{Fn}
    f = F\cdot n = (\tau\times (\sigma b))\cdot n = (n\times \tau)\cdot(\sigma b)
        = \nu\cdot(\sigma b) = (\sigma\nu) \cdot b.
\end{equation}

The Peach-Koehler formula, however, has the following mathematical problem. If
we treat the dislocation as a one-dimensional curve $\Gamma(t)$, then the
stress field $\sigma$ mathematically has a singularity along the dislocation
curve as we see in Sect.~\ref{subsec:2} and there is no pointwise value of
$\sigma$ on $\Gamma(t)$. Furthermore, the stress field must have infinite
energy under a most naive setting of the problem if the dislocation curve is a
mathematically sharp one-dimensional object. These mathematical singularities
seem to be physically regularized due to the existence of a minimum scale given
by the size of atoms and lattice spacing. In physics, this problem is often
solved by introducing the notion of dislocation core which is a tubular region
around the dislocation line of thickness of about $5|b|$ (see
\cite{Bulatov06}).

In this paper, we consider a regularized mathematical model of the motion of
dislocation curves by means of phase field modeling. The model is derived as
a gradient flow of an energy in a mathematically systematic way, and we show
that it naturally includes the Peach-Koehler force in Sect.~\ref{sec:2}. A
simplified 2D-1D model is also derived in Sect.~\ref{sec:3} and its numerical
examples are presented in Sect.~\ref{sec:4}.

\section{An energy-based approach to modeling of dislocation dynamics}
\label{sec:2}

In this section, we study displacement field in an elastic body with a
dislocation, and observe that the total elastic energy is infinite without any
regularization. A phase field model is proposed by introducing a regularized
energy in Sect.~\ref{regularization}.

\subsection{Elastic energy with dislocation}
\label{subsec:2}

We start from a homogeneous elastic body without dislocation. We denote it by
$\Omega$ which is a bounded Lipschitz domain in $\R^3$. The position vector in
$\bar{\Omega}$ is denoted by $x=(x_1,x_2,x_3)\uT \in\bar{\Omega} \subset\R^3$,
where $\uT$ denotes the transpose of a vector or matrix. All vectors are assumed
to be column vectors in this paper. We use the following notation: $\partial_j
:= \frac{\partial}{\partial x_j}$, $\nabla = (\partial_1, \partial_2,
\partial_3)\uT$, $\nabla\uT u = (\partial_j u_i)\in \Tensor$ and $\nabla u\uT =
(\nabla\uT u)\uT$ for $u(x)\in\R^3$. We often use Einstein's summation
convention for the space variables. For matrices $A = (a_{ij})$, $B = (b_{ij})\in
\R^{3\times 3}$, their inner product is denoted by $A:B := a_{ij}b_{ij}$.

Small deformation of the elastic body is described by a displacement field
$u(x) = \big(u_1(x),u_2(x),u_3(x)\big)\uT \in\R^3$, the symmetric strain
tensor $e[u](x) = (e_{ij}[u](x))$:
\begin{equation}\label{s-d-relation}
    e[u] := \frac{1}{2}\big(\nabla\uT u+\nabla u\uT\big) \in\SymTensor, \quad\text{i.e.},\quad
        e_{ij}[u] = \frac{1}{2}\big(\partial_ju_i+\partial_iu_j\big),
\end{equation}
and the stress tensor $\sigma[u](x) = \big(\sigma_{ij}(x)\big) \in\SymTensor$:
\[
    \sigma[u](x) = C(x)e[u](x), \quad \text{i.e.,} \quad
        \sigma_{ij}(x) = c_{ijkl}(x) e_{kl}(x),
\]
where $C(x) = \big(c_{ijkl}(x)\big) \in\R^{3\times 3\times 3\times 3}$ is the
(anisotropic) elasticity tensor with the symmetries $c_{ijkl} = c_{klij} =
c_{jikl}$, $i,j,k,l\in \{ 1,~2,~3\}$. It should satisfy the positivity
condition:
\begin{equation}\label{pcondition}
    ^\exists c_*>0 \quad \text{s.t.} \quad c_{ijkl}(x) \xi_{ij} \xi_{kl} \geq c_*|\xi|^2,
        \qquad ^\forall x\in\Omega_0,\ ^\forall \xi\in\SymTensor,
\end{equation}
where $|\xi | := \sqrt{\xi_{ij}\,\xi_{ij}}$. It depends on the elastic property
of the material $\Omega$ and is supposed to be given. If the material is
homogeneous, the elasticity tensor should be constant $C(x) \equiv C$. From the
strain-displacement relation \eqref{s-d-relation}, we write $\sigma [u] := C \, e[u]$.

The displacement field $u(x)$ is obtained by the following linear second order
elliptic boundary value problem:
\begin{equation}\label{elasticity}
\left\{
    \begin{aligned}
        -\partial_j \sigma_{ij}[u] &= f_i(x), &i=1,2,3,\ x&\in\Omega,\\
        u &= g(x), &x&\in\GammaD,\\
       \sigma[u] \nu &= \fN(x), &x&\in\GammaN.
  \end{aligned}
\right.
\end{equation}
where $f(x) =\big(f_1(x), f_2(x), f_3(x)\big)\uT \in\R^3$ is a given body
force, and the first equation represents the equilibrium equations of force.
The boundary $\partial \Omega$ is divided into two parts as
\[
    \partial\Omega =\GammaD\cup \GammaN, \qquad \GammaD\cap\GammaN = \emptyset,
        \qquad |\GammaD|>0,
\]
where $\GammaD$ is an open non-empty portion of $\Gamma$. The two dimensional
area of $\Gamma_D$ is denoted by $|\GammaD|$. On the other hand, $\GammaN$ can
be empty. The outward unit normal vector to $\partial\Omega$ at
$x\in\partial\Omega$ is denoted by $\nu(x)\in\R^3$. The displacement on
$\GammaD$ is given by $g\in H^{\1/2}(\GammaD; \R^3)$, and the surface outer
force on $\GammaN$ is given by $\fN \in L^2(\GammaD;\R^3)$, as prescribed
boundary values.

Under suitable regularity conditions, it is well-known that a unique solution
to the problem \eqref{elasticity} is given by the unique minimizer of a total
elastic energy
\begin{equation}\label{ele}
    E_1(u) := \int_\Omega \left\{ w[u] - f\cdot u \right\} \dd x 
       - \int_{\GammaN} \fN \cdot u \dd s,
\end{equation}
with the boundary condition $u = g$ on $\GammaD$, where $w[u]$ is the strain
energy density defined as $w[u]:=\1/2 \sigma [u]:e[u]$.

In case that the elastic material is isotropic and homogeneous, the elasticity
tensor has the form
\[
    c_{ijkl} = \lambda \delta_{ij}\delta_{kl} + \mu(\delta_{ik}\delta_{jl} + \delta_{il}\delta_{jk}),
\]
where $\lambda\ge 0$ and $\mu >0$ are called the Lam\'{e} constants. Since
\[
    c_{ijkl}\xi_{ij}\xi_{kl} = \lambda (\mathrm{tr}\,\xi)^2 + 2\mu|\xi|^2,
        \qquad \xi \in\SymTensor,
\]
the condition \eqref{pcondition} is satisfied with $c_* = 2\mu$. The stress
tensor and the strain energy density become
\[
    \sigma[u] = \lambda (\Div u) I + 2\mu e[u], \qquad
        w[u] = \1/2 \lambda (\Div u)^2 + \mu |e[u]|^2,
\]
where $I$ denotes the unit tensor. The equilibrium equations of force are often
called the Navier or Navier-Cauchy equations and take the following form:
\begin{equation}\label{navier}
    -\mu \Delta u - (\lambda +\mu )\nabla (\Div u)= f,
        \qquad \text{in }\Omega.
\end{equation}

Let us now consider the case that $\Omega$ contains a dislocation curve
$\Gamma(t)$. We assume that $\Gamma(t)$ is a closed plane curve without
self-intersections in a fixed crystallographic plane $\tilde{S}\subset \R^3$
and define $S := \tilde{S}\cap \Omega$, which is called the slip plane.  We
suppose that $S$ is connected and open in $\tilde{S}$.  The Burgers vector of
$\Gamma(t)$ is denoted by $b\in \R^3\setminus\{\bm{0}\}$ which is a fixed
vector tangential to $S$. This is a so-called mixed dislocation which contains
both edge and screw dislocations at the parts of $\Gamma(t)$ where $b \perp
\Gamma(t)$ and $b \parallel \Gamma(t)$, respectively.

Choosing a suitable orthogonal coordinate system, without loss of generality,
we suppose that $\tilde{S} = \{x=(x_1,x_2,x_3)\uT;\ x_3=0\}$, $b =
(b_1,0,0)\uT$, and $b_1>0$.  A typical example of $\Omega$ and $S$ is shown in
Fig.~\ref{fig:1} and Fig.~\ref{fig:2}, where 
\[
    \Omega =S_0\times (-H,H),
    \qquad
    \GammaD =S_0\times \{-H,H\},
    \qquad
    \GammaN =\partial S_0 \times [-H,H].
\]
We often identify the slip plane $S = S_0\times \{0\}$ with $S_0\subset \R^2$,
if no confusion occurs.  The coordinate in $S\cong S_0\subset \R^2$ is denoted
by $x'=(x_1,x_2)\uT$.  We consider $\Gamma(t)$ and $S$ as subsets of $\R^2$,
and a two dimensional domain enclosed by $\Gamma(t)$ in $S$ is denoted by
$D(t)\subset S\subset \R^2$.

We define $\Omega^\pm := \{(x_1,x_2,x_3)\uT\in\Omega;\ \pm x_3 > 0\}$ and
denote the outward unit normal vectors on $\partial \Omega^\pm$ by $\nu^\pm$,
respectively. It is considered that the displacement field $u$ is discontinuous
across the slip plane $S$. The traces of $u$ to $S$ from $\Omega^+$ and
$\Omega^-$ are denoted by $u^+$ and $u^-$, respectively, and normal tractions
$\sigma [u] \nu$ from $\Omega^+$ and $\Omega^-$ on $S$ are denoted by
$\sigma^+\nu$ and $\sigma^-\nu$, respectively, where we define $\nu =
(0,0,1)\uT$ on $S$.  We also denote the outward unit normal vector on $\partial
\Omega$ by $\nu$. We remark that $\nu = \nu^\pm$ on $\partial \Omega\cap
\partial \Omega^\pm$ and that $\nu =\mp \nu^\pm$ on $S$. The gaps of the
displacement and the traction across the slip plane $S$ are denoted by
$\Jump{u} := u^+-u^-$ and $\Jump{\sigma} \nu := \sigma^+\nu -\sigma^-\nu$,
respectively.

For a fixed time $t$, it is naturally expected that the displacement field
$u(x)$ in $\Omega$ satisfies the following boundary value problem in a naive
setting:
\begin{equation}\label{jumpprob}
\left\{
\begin{aligned}
    -\partial_i \sigma_{ij}[u] &= f_j(x),\quad j = 1,2,3,& x &\in\Omega^+\cup\Omega^-,\\
    u &= g(x),& x&\in\GammaD,\\
    \sigma[u]\nu &= \fN(x),& x &\in\GammaN,\\
    \Jump{u} &=
        \begin{cases}
            b, & \quad \text{on}\ D(t),\\
            0, & \quad \text{on}\ S\setminus D(t),
        \end{cases}\\
    \Jump{\sigma[u]}\nu &= 0, \qquad \text{on}\ S.
\end{aligned}
\right.
\end{equation}
This problem, however, has no finite energy solution as seen below.  Let us
suppose that $u$ is a finite energy solution, i.e., $u|_{\Omega^\pm}\in
H^1(\Omega^\pm;\R^3)$.  Then its traces $u^\pm$ and $\Jump{u}$ should belong to
$H^{1/2}(S;\R^3)$ but it contradicts the fourth condition of
\eqref{jumpprob}.  In the next section, we consider a weak formulation of a
jump problem in general form and study the condition that its solution has a
finite energy.

\subsection{Weak formulations for jump problems}
\label{jp}

We study a weak formulation for a linear elasticity problem with a jump
condition across an interface like \eqref{jumpprob} in a general setting.

Under the same condition of \eqref{jumpprob}, for a $\R^3$-valued function
$\Psi$ on $S$, we consider the following problem
\begin{equation}\label{gjp}
\left\{
\begin{aligned}
    -\partial_j \sigma_{ij}[u] &= f_i(x), \quad i=1, 2, 3,& x &\in\Omega^+\cup\Omega^-,\\
    u &= g(x),& x &\in\GammaD,\\
    \sigma[u]\nu &= \fN(x),& x &\in\GammaN,\\
    \Jump{u} &= \Psi,& &\text{on } S,\\
    \Jump{\sigma[u]} \nu &= 0,& &\text{on } S.
\end{aligned}
\right.
\end{equation}

We suppose that $\Psi\in H^{\1/2}(S;\R^3)$ and $g\in H^{\1/2}(\GammaD;\R^3)$.
For the problem \eqref{gjp}, we define function spaces:
\begin{align*}
X &:= \big\{ w\in L^2(\Omega;\R^3);\ w|_{\Omega^\pm}\in H^1(\Omega^\pm;\R^3) \big\},\\
X(\Psi) &:= \big\{ w\in X;\ \Jump{w} = \Psi \big\},\\
V(g,\Psi) &:= \big\{ w\in X(\Psi);\ w=g \text{ on } \GammaD \big\},\\
V &:= V(0,0),
\end{align*}
where $\Jump{w} = w^+-w^-$ on $S$ for $w\in X$. The space $X$ is a Hilbert
space with the following norm and a corresponding inner product:
\[
    \|w\|_X := \left(\|w\|_{H^1(\Omega^+;\R^3)}^2 + \|w\|_{H^1(\Omega^-;\R^3)}^2\right)^\1/2.
\]
We can identify $X(0) = H^1(\Omega;\R^3)$. We remark that $X$ is a space of
finite energy displacements with a gap across $S$. In the following lines, we
suppose that $V(0,\Psi)\neq\emptyset$.

If $u|_{\Omega^\pm}\in H^2(\Omega^\pm;\R^3)$, and $u$ satisfies the equations
of \eqref{gjp}, $u$ is called a strong solution to \eqref{gjp}, where the
boundary conditions are considered in the sense of the trace operator.  We also
define a weak solution as follows.
\begin{Prob}\label{weakjp}
Find $u \in V(g,\Psi)$ such that
\begin{equation}\label{ws}
    \int_{\Omega\setminus S} \sigma[u]:e[w] \dd x
        - \int_\Omega f\cdot w \dd x
        - \int_{\GammaN}\fN\cdot w \dd s = 0,\quad \text{for all}\quad w\in V.
    \end{equation}
\end{Prob}
A solution of Problem~\ref{weakjp} is called a weak solution to \eqref{gjp}.
In particular, it is not difficult to show that a strong solution to \eqref{gjp}
is a weak solution, by a standard computation with integration by parts.
More precisely, we have the following theorem.
\begin{Th}\label{equi1}
    Suppose that $f\in L^2(\Omega;\R^3)$, $\fN\in L^2(\GammaN,\R^3)$, $g\in
    H^{\1/2}(\GammaD;\R^3)$ and $\Psi\in H^{\1/2}(S;\R^3)$.  Then, $u$ is a
    strong solution to \eqref{gjp} if and only if $u$ is a weak solution and
    $u|_{\Omega^+}$ and $u|_{\Omega^-}$ belong to $H^2(\Omega^+;\R^3)$ and
    $H^2(\Omega^-;\R^3)$, respectively.
\end{Th}
The unique existence of the weak solution is guaranteed as follows.
\begin{Th}\label{equi2}
    Under the same conditions of Theorem~\ref{equi1}, there exists a unique
    solution $u$ to Problem~\ref{weakjp} and $u$ is given as a unique minimizer of
    $E_1(u)$ in $V(g,\Psi)$. 
\end{Th}
We omit the proofs of these theorems here and postpone them until our
forthcoming paper. Here we just admit that the weak solution uniquely exists
and let our argument proceed.

%


On the other hand, if $\Psi$ does not belong to $H^{1/2}(S;\R^3)$ as in the
dislocation model \eqref{jumpprob}, there is no finite energy solution since
$\Psi$ cannot be expressed by $\Jump{u}$ for any $u\in X$. According to the
model \eqref{jumpprob}, the displacement field $u$ has to have a singularity
along the dislocation curve with infinite elastic energy, and of course this is
not a realistic solution. In a microscopic description of a crystal lattice, it
is considered that this singularity is somehow regularized by the existence of
a minimum length corresponding to the height of an atom. In the next section,
we introduce a regularization in terms of phase field approach.


\subsection{Regularization of energy}\label{regularization}

We consider a regularization for the curvature flow model \eqref{cf} by means
of an Allen-Cahn-type bistable potential energy.  Let $\varphi$ be a phase
field which is a smooth scalar-valued function defined on $S\times [0,T)$ and
the value is approximately $1$ in $D(t)$ and $0$ in $S\setminus
\overline{D(t)}$.  For positive parameters $\vep >0$ and $\beta >0$, we define
an interface energy as
\begin{equation}\label{ie1}
    E_0(\varphi) := \int_S \left(\frac{\vep}{2}|\nabla' \varphi|^2
        + \beta W(\varphi)\right) \dd x',
\end{equation}
where $W(s) = s^2(1-s)^2$. As its gradient flow in $L^2(S)$, we have the following
Allen-Cahn equation
\begin{equation}\label{ac}
    \alpha \frac{\partial \varphi}{\partial t} = \vep \Delta' \varphi -\beta W'(\varphi),
        \qquad \text{in}\ S\times (0,T),
\end{equation}
where $\nabla'$ and $\Delta'$ represent the gradient and the Laplacian with
respect to $x'=(x_1,x_2)\uT$. It is known that the interface energy \eqref{ie1}
and the Allen-Cahn equation correspond to the length of the curve $\Gamma(t)$
and a motion by line tension, respectively, under a suitable scaling
\cite{Fife88}.

Using this regularized interface energy $E_0(\varphi)$ together with the
elastic energy $E_1(u)$ defined by \eqref{ele}, we consider the following total
energy:
\[
    E(\varphi) := E_0(\varphi) + E_1(u_\varphi), \qquad \varphi\in H^1(S),
\]
where $u_\varphi$ is a unique solution of \eqref{ws} with $u_\varphi\in
V(g,b\varphi)$.

Let us derive an $L^2$ gradient flow of the energy $E(\varphi)$. For a smooth
scalar function $\psi$ defined on $S$, we consider a first variation:
\begin{align}\label{fv}
    \frac{d}{dr} E(\varphi+r \psi)|_{r=0}
        = \frac{d}{dr} E_0(\varphi+r \psi)|_{r=0} + \frac{d}{dr} E_1(u_{\varphi+r \psi})|_{r=0}.
\end{align}
The first term is formally given as
\[
    \frac{d}{dr} E_0(\varphi+r \psi)|_{r=0}
        = \int_S \big(-\vep \Delta' \varphi + \beta W'(\varphi) \big)\psi \dd x'
            + \int_{\partial S}\frac{\partial \varphi}{\partial \nu'}\psi \dd s',
\]
where $\nu'\in \R^2$ denotes the outward unit normal vector on $\partial
S\cong\partial S_0$.  For the second term of \eqref{fv}, it is easy to show
that $u_{\varphi +r  \psi} = u_\varphi+r  u_* \in V(g,b(\varphi+r \psi))$ where
$u_*$ is a unique solution of the following weak form:
\begin{equation}\label{ws0}
    u_*\in V(0,b\psi), \qquad
    \int_{\Omega\setminus S} \sigma[u_*]:e[w] \dd x = 0, \quad \text{for all}\quad w\in V.
\end{equation}
Then, under suitable regularity conditions, we obtain
\begin{align*}
    \frac{d}{dr }E_1(u_{\varphi+r \psi})|_{r=0}
        &= \frac{d}{dr }E_1(u_\varphi +r u_* )|_{r=0}\\
        &= \int_{\Omega^+\cup\Omega^-} \big\{ \sigma[u_\varphi]:e[u_* ] - f\cdot u_* \big\} \dd x 
            - \int_{\GammaN} \fN\cdot u_*  \dd s\\
        &= \int_{\partial \Omega^+}(\sigma[u_\varphi]\nu^+)\cdot u_* \dd s
            + \int_{\partial \Omega^-}(\sigma[u_\varphi]\nu^-)\cdot u_* \dd s\\
        &\qquad - \int_{\Omega^+\cup\Omega^-} \big\{ \Div\sigma[u_\varphi]+f\big\}\cdot u_*  \dd x
            - \int_{\GammaN} \fN\cdot u_*  \dd s\\
        &= \int_{S} \big\{(-\sigma^+[u_\varphi]\nu)\cdot u_* ^+ +(\sigma^-[u_\varphi]\nu)\cdot u_* ^-\big\}\dd x' \\
        &\qquad + \int_{\GammaN} (\sigma[u_\varphi]\nu - \fN)\cdot u_*  \dd s\\
        &= -\int_{S} (\sigma[u_\varphi]\nu)\cdot \Jump{u_* } \dd x'\\
        &= -\int_{S} (\sigma[u_\varphi]\nu)\cdot b \psi \dd x'.
\end{align*}
Hence, using the homogeneous Neumann boundary condition $\frac{\partial
\varphi}{\partial \nu'}=0$ on $\partial S$, we obtain
\[
    \frac{d}{dr} E(\varphi+r \psi)|_{r=0}
    = \int_S \big( -\vep \Delta' \varphi +\beta W'(\varphi) -(\sigma[u_\varphi]\nu)\cdot b \big)\psi \dd x'.
\]
Similarly to the Allen-Cahn equation \eqref{ac}, we derive a gradient flow of
the energy $E(\varphi)$. Therefore, we propose the following phase field model
for dislocation dynamics
\begin{align}\label{pfm}
\left\{
\begin{aligned}
    \alpha \frac{\partial \varphi}{\partial t} &= \vep \Delta' \varphi -\beta W'(\varphi) +(\sigma[u_\varphi]\nu)\cdot b,&
        \text{in} &\ S\times (0,T),\\
    \frac{\partial \varphi}{\partial \nu'} &= 0,& \text{on} &\ \partial S\times (0,T),\\
    \varphi(\cdot,0) &= \varphi_0, &\text{on} &\ S,
\end{aligned}
\right.
\end{align}
where $\varphi_0(x')$ is a suitable initial value for $\varphi$.  The
additional force term $(\sigma[u_\varphi]\nu)\cdot b$ appearing above is
nothing but the Peach-Koehler force \eqref{Fn} acting on the dislocation line.
The model \eqref{pfm} naturally contains the notion of the dislocation core and
the Peach-Koehler force. The gradient flow structure
\begin{equation}\label{gfs}
    \frac{d}{dt}E(\varphi(\cdot,t)) = -\alpha \int_{S}
    \left|\frac{\partial \varphi}{\partial t}(x',t)\right|^2\,dx' \le 0,
\end{equation}
behind the dynamics of dislocation curve has been revealed in the derivation of
the model as above.

\begin{figure}
\includegraphics[width=1.0\linewidth]{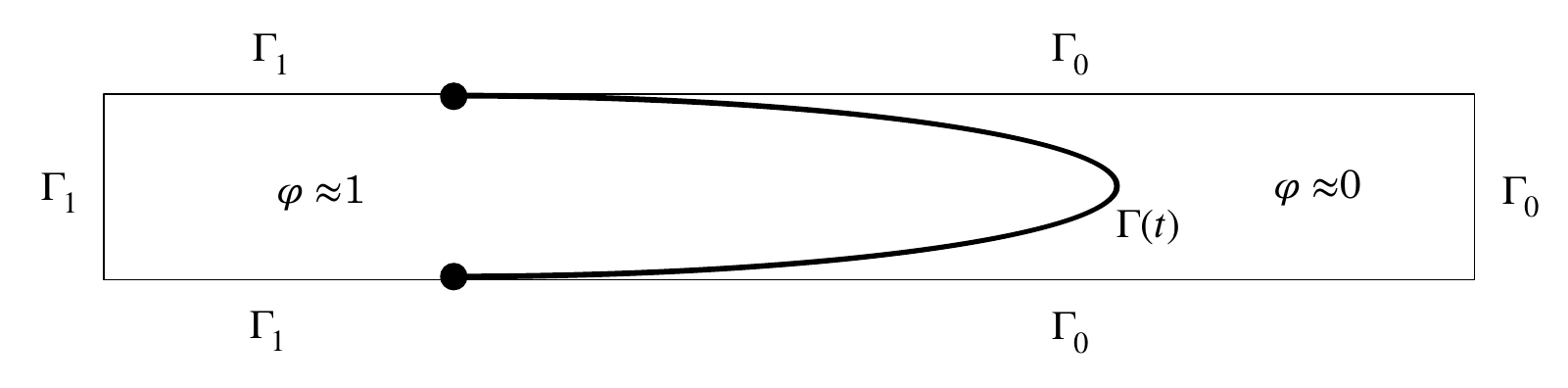}
\caption{Another possible configuration of the dislocation curve with 
fixed end points in our model is shown. The Dirichlet boundary condition for $\varphi$ is imposed as 
$\varphi |_{\Gamma_0}\approx 0$ and 
$\varphi |_{\Gamma_1}\approx 1$.}
\label{fig:with_endpoints}
\end{figure}

\begin{Rem}
In this paper, we treat only the case where the dislocation is a closed
plane curve without any end point in the slip plane.
On the other hand, it is often observed in real materials that two end points
are fixed at some defect of the crystal structure.
It is also known that the dislocation curve sometimes can change its slip plane.
For example, some numerical simulations of dislocations with end points which
change their slip planes are shown in \cite{Paus12}.
The model presented in this paper can also be applied in the case where the end
points are fixed on the boundary of the slip plane, if we impose the Dirichlet
boundary condition of the phase field variable $\varphi$ as illustrated in
Figure~\ref{fig:with_endpoints}.
However, the treatment of the change of the slip plane as shown in
\cite{Paus12} in our model seems challenging but difficult at present.
\end{Rem}

\section{A simplified 2D-1D model}\label{sec:3}

In this section, we derive a 2D-1D coupled phase field model for dynamics of a
straight screw dislocation line in the same manner as the 3D-2D model of the
previous section.

\begin{figure}[b]
\includegraphics[width=0.6\linewidth]{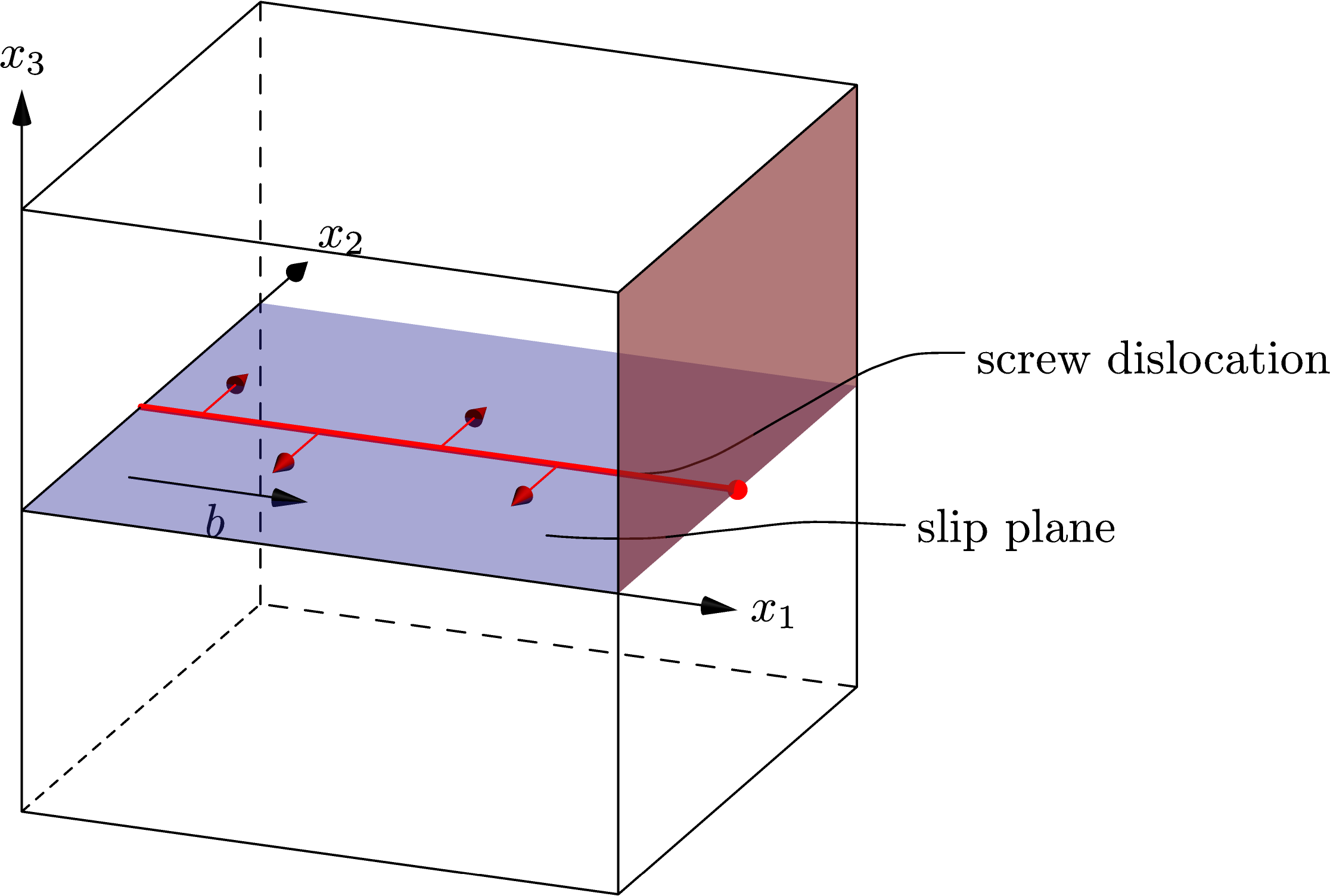}
\raisebox{1cm}{\includegraphics[width=0.38\linewidth]{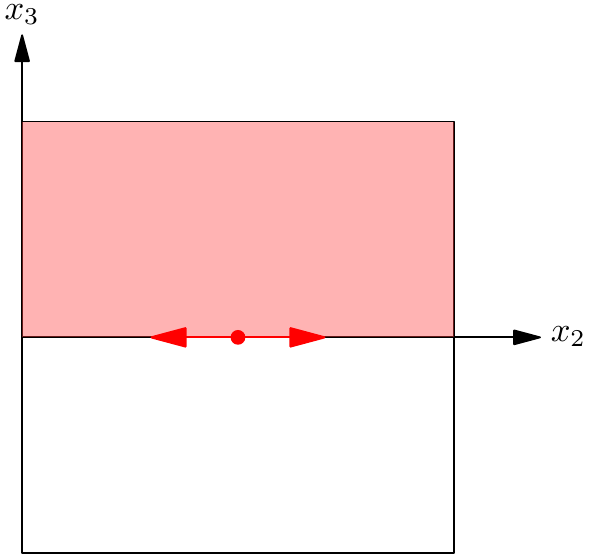}}
\caption{A straight screw dislocation in a slip plane and its two dimensional
projection.}
\label{fig:screw}
\end{figure}

As shown in Fig.~\ref{fig:screw}, we consider a rectangular parallelepiped
elastic body and a slip plane in $x_1x_2$ plane.  A straight screw dislocation
with the Burgers vector $(b,0,0)\uT$ is assumed to be moving in the slip plane.
In the following lines, we denote the coordinate $(x_2,x_3)$ also by $(x,y)$
for simplicity.  For the displacement field, we suppose the so-called
anti-plane displacement and an odd symmetry with respect to the slip plane
$y=0$, i.e., the displacement has the form $(u(x,y), 0,0)\uT$ and
$u(x,-y)=-u(x,y)$.  Then the equilibrium equation \eqref{navier} becomes $-\mu
\Delta u =0$, where no body force is supposed.

In the following sections, we set a rectangular two-dimensional domain
$\Omega:=(-L,L)\times (0,H)$ in $xy$ plane. We define the bottom boundary
$S:=(-L,L)\times \{0\}$, the lateral boundary $\Gamma_0:=\{-L,L\}\times [0,H]$, and
top boundary $\Gamma_1:=(-L,L)\times \{H\}$. Then $\partial \Omega = S\cup
\Gamma_0\cup\Gamma_1$ holds. The Laplacian with respect to $(x,y)$ is denoted
by $\Delta =\frac{\partial^2}{\partial x^2} +\frac{\partial^2}{\partial y^2}$.

Let $u(x,y,t)\in\R$ be an anti-plane displacement at time $t$ and position
$(x,y)\in\Omega$.  Similarly to Sect.~\ref{regularization}, we introduce a
phase field variable $\varphi\in H^1(-L,L)$ and define the following total
energy.  We suppose that an anti-plane boundary traction $g(x)$ is given on the
top boundary $\Gamma_1$. We set $F(s)=s^2(1-s)^2$ and define
\begin{align*}
    E(\varphi) &:= E_0(\varphi) + E_1(u_\varphi),\\
    E_0(\varphi) &:= \int_{-L}^L \left(\frac{\vep}{2}|\varphi_x|^2
        + \beta F(\varphi)\right) \dd x,\\
    E_1(u) &:= \mu\int\hspace{-5pt}\int_\Omega |\nabla u|^2 \dd x \dd y
        -2\int_{\Gamma_1}gu \dd x,
\end{align*}
where $u_\varphi \in H^1(\Omega)$ is given as the unique weak solution to
\begin{equation}\label{lp}
\left\{
\begin{aligned}
    \Delta u &= 0,& \text{ in }&\Omega,\\
    u_x &= 0, &\text{ on }&\Gamma_0,\\
    \mu u_y &= g(x),& \text{ on }&\Gamma_1,\\
    u &= \frac{b}{2}\varphi(x),& \text{ on }&S.
\end{aligned}
\right.
\end{equation}

Let us derive an $L^2$-gradient flow of the energy $E(\varphi)$. Under the
Neumann boundary condition $\varphi_x(-L) = \varphi_x(L) = 0$,
for a smooth function $\psi$ defined on $S$, we consider a first variation
\begin{equation}\label{fv2}
    \frac{d}{dr} E(\varphi+r \psi)|_{r=0}
    =
    \frac{d}{dr} E_0(\varphi+r \psi)|_{r=0}
    +
    \frac{d}{dr} E_1(u_{\varphi+r \psi})|_{r=0}.
\end{equation}
The first term is formally given as
\[
    \frac{d}{dr} E_0(\varphi+r \psi)|_{r=0}
    =
    \int_{-L}^L \left(-\vep \varphi_{xx} + \beta F'(\varphi) \right)\psi \dd x.
\]

It is easy to show that $u_{\varphi + r \psi} = u_\varphi + \frac{1}{2} r b w
\in V(g,\varphi+r \psi)$, where $w\in H^1(\Omega)$ is the unique weak solution of 
\begin{equation*}
\left\{
\begin{aligned}
    \Delta u &= 0,& \text{ in }&\Omega,\\
    u_x &= 0, &\text{ on }&\Gamma_0,\\
    u_y &= 0,& \text{ on }&\Gamma_1,\\
    u &= \psi(x),& \text{ on }&S.
\end{aligned}
\right.
\end{equation*}
Then the second term of \eqref{fv2} becomes
\begin{align*}
    \frac{d}{dr }E_1(u_{\varphi+r  \psi})|_{r=0}
    &= \frac{d}{dr }E_1(u_\varphi +\frac{1}{2}r  bw)|_{r=0}\\
    &= \mu b \iint_{\Omega} \nabla u_\varphi \cdot \nabla w \dd x \dd y
        - b \int_{\Gamma_1} g w\dd x\\
    &= \mu b\int_{\partial \Omega}\frac{\partial u_\varphi}{\partial \nu} w \dd s
        -\mu b \iint_{\Omega} (\Delta u_\varphi) w \dd x \dd y
        - b \int_{\Gamma_1} g w \dd x\\
    &= -\mu b\int_S \frac{\partial u_\varphi}{\partial y} \psi \dd x.
\end{align*}
Hence, we obtain
\[
    \frac{d}{dr} E(\varphi+s \psi)|_{r=0}
    = \int_{-L}^L \left(-\vep \varphi_{xx} +\beta F'(\varphi)
        -\mu b \frac{\partial u_\varphi}{\partial y}(x,0) \right)\psi \dd x.
\]
Similarly to the Allen-Cahn equation \eqref{ac}, we derive a gradient flow of
the energy $E(\varphi)$ as
\begin{equation}\label{phieq}
    \alpha \varphi_t = \vep \varphi_{xx} -\beta F'(\varphi)
        + b\mu \frac{\partial u_\varphi}{\partial y}(x,0), \qquad x\in (-L,L),\ t>0.
\end{equation}
The additional force term $b\mu \frac{\partial u_\varphi}{\partial y}(x,0)$ is
formally corresponding to the Peach-Koehler force \eqref{Fn} acting on the
screw dislocation. We substitute the relation
\[
    \varphi(x,t) = \frac{2}{b}u_\varphi(x,0,t),\qquad x\in [-L,L],~t\ge 0,
\]
into \eqref{phieq} and set
\[
    W(s) := \frac{4\beta}{b^2}s^2\left(\frac{b}{2}-s\right)^2,
    \qquad
    \gamma := \frac{\mu b^2}{2}.
\]
The boundary load $g(x)$ is also assumed to depend on time $t$. Then, we obtain
the following 2D-1D dislocation model:
\begin{equation}\label{themodel}
\left\{
\begin{aligned}
    \Delta u &= 0,&  (x,y)&\in\Omega,\ t\in [0,T],\\
    u_x &= 0,& (x,y)&\in\Gamma_0,\ t\in [0,T],\\
    u_y &= g(x,t),& (x,y)&\in\Gamma_1,\ t\in [0,T],\\
    \alpha u_t &= \vep u_{xx} - W'(u) + \gamma u_y,& (x,y)&\in S,\ t\in(0,T],\\
    u &= u_0(x),& (x,y)&\in S,\ t=0.
\end{aligned}
\right.
\end{equation}

\section{Numerical results}\label{sec:4}

In this section, we give some numerical examples for an approximation of
\eqref{themodel}.

We discretize \eqref{themodel} by employing standard finite differences to
approximate spatial derivatives. A rectilinear grid is obtained by dividing
$\Omega$ into $N_x\times N_y$ rectangular cells and a solution is sought at the
grid nodes. This spatial discretization results in a non-linear system of ODEs
in time that we solve by means of a fully implicit, variable-step solver
\cite{CVode}. 

The initial condition in both examples below is a shifted and scaled Heaviside
function $u_0(x) = b H(x - x_0)$ where $x_0$ is the initial position of the
step. Other parameter values used in the numerical examples are summarized in
Table~\ref{tab:params}.

In both figures below the horizontal axis corresponds to the spatial $x$
variable while the vertical axis corresponds to the time. We plot the graph of
$u$ at $S$ at 50 time levels distributed uniformly in $[0,T]$, the initial
condition $u_0$ being the lowest graph with the time increasing upwards. The
graph of $u$ in the whole of $\Omega$ is not shown.

\newcolumntype{Y}{>{\centering\arraybackslash}X}
\begin{table}
    \centering
    \begin{tabularx}{\linewidth}{*{10}{Y}}
        \toprule
        $T$ & $L$ & $H$ & $N_x$ & $N_y$ & $b$ & $\varepsilon$ & $\beta$ & $\mu$ & $\alpha$ \\
        \midrule
        4 & 2 & 2 & 256 & 128 & 0.06 & 0.04 & 10 & 10 & 0.01 \\
        \bottomrule
    \end{tabularx}
    \caption{Parameter values used in numerical examples for the model \eqref{themodel}.}
    \label{tab:params}
\end{table}

\begin{figure}
\centering
\includegraphics[width=0.8\linewidth]{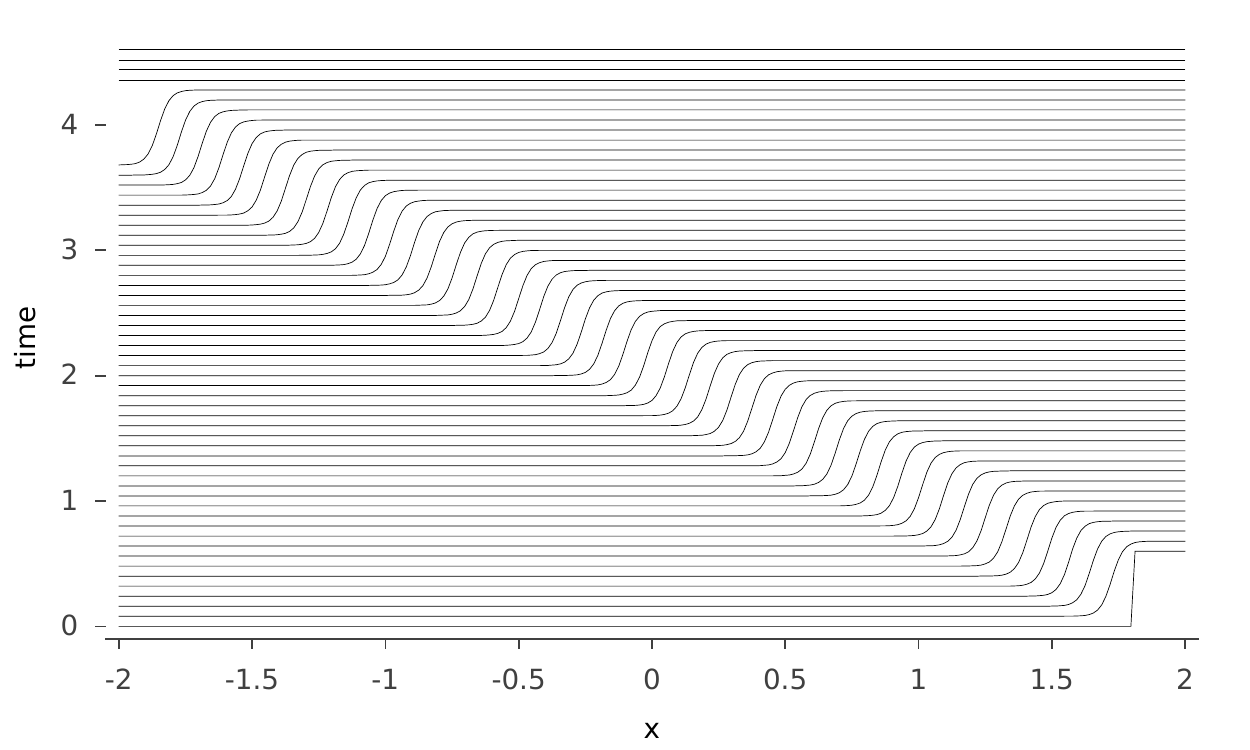}
\caption{Numerical example of a screw dislocation moving under a constant loading.}
\label{fig:3}
\end{figure}

In Fig.~\ref{fig:3} we present a numerical example of a screw dislocation that
moves in a stress field caused by a constant anti-plane traction imposed on the
top boundary $\Gamma_1$. We set $g\equiv 0.5$ and place the dislocation at the
initial location $x_0 = 1.8$. The dislocation moves to the left at an almost
constant velocity until it reaches the crystal surface at $x=-2$ where it
annihilates.

\begin{figure}
\centering
\includegraphics[width=0.8\linewidth]{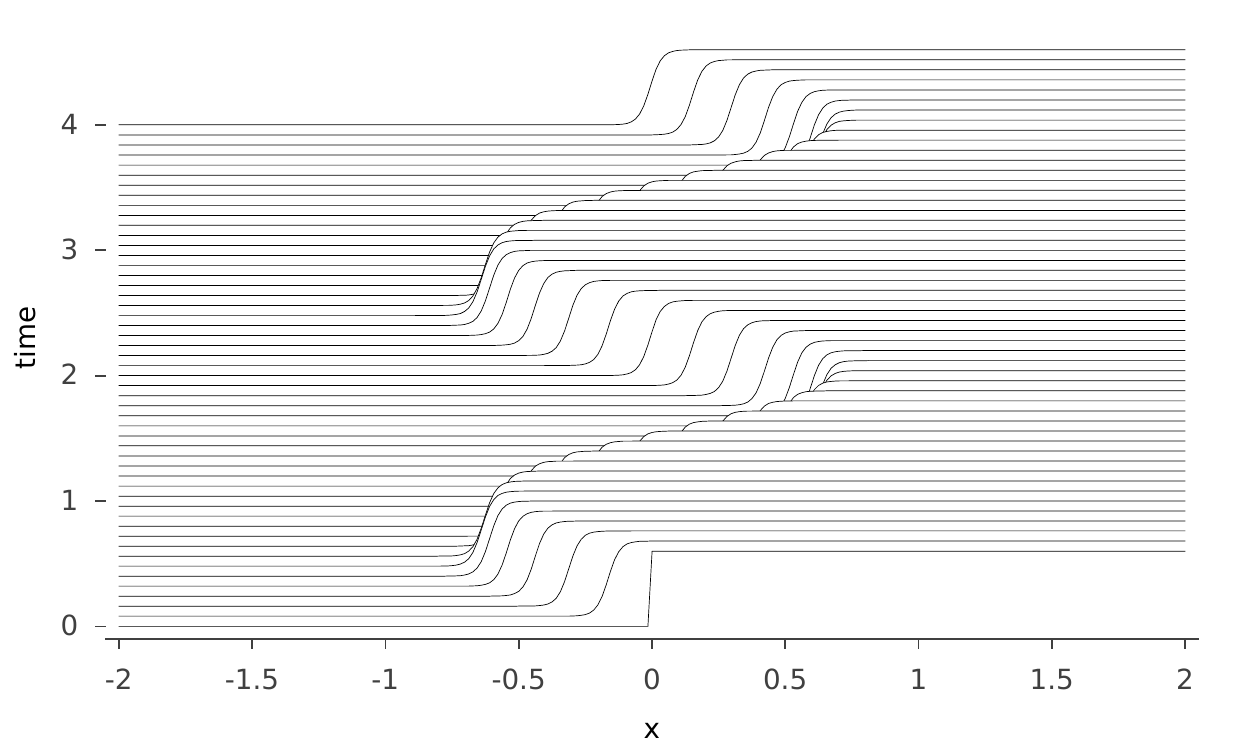}
\caption{Numerical example of a screw dislocation moving under a periodic loading.}
\label{fig:4}
\end{figure}

In Fig.~\ref{fig:4} we consider the behaviour of a screw dislocation under a
periodic loading. We place the dislocation at $x_0 = 0$ and set $g(x,t) =
\cos(0.5t)$. The oscillating stress field causes the dislocation to move
periodically around its initial location.

\section{Summary and Conclusions}\label{conclusions}

We proposed a 3D-2D phase field model \eqref{pfm} for dynamics of a mixed
dislocation curve and a 2D-1D model \eqref{themodel} for a straight screw
dislocation. They are both derived as gradient flows of regularized total
energies and naturally include the Peach-Koehler force and the notion of the
dislocation core. Some numerical examples of the 2D-1D model are given in
Sect.~\ref{sec:4}. The revealed gradient flow structure \eqref{gfs} is expected
to be useful for further mathematical and numerical analysis.

One of interesting questions about this model is the relation between the layer
width of $\varphi$ (i.e., the radius of the dislocation core) and the
parameters $\varepsilon$, $\beta$ and $|b|$. This is still open but seems to be
important not only in the sense of modeling but also in numerical simulations
for choosing a proper mesh size.

\end{document}